

\baselineskip=14pt
\parskip=10pt

\magnification=\magstephalf

\def\1{{\overline{1}}}
\def\2{{\overline{2}}}
\parindent=0pt
\overfullrule=0in

\def\frac#1#2{{#1 \over #2}}
\centerline
{\bf A Detailed Analysis of Quicksort Running Time} 
\bigskip
\centerline
{\it Shalosh B. EKHAD and Doron ZEILBERGER}
\bigskip

{\bf Quicksort: Sir Tony Hoare's Seminal Invention}

In 1905, 25-year-old Albert Einstein revolutionized physics, while in 1931, 25-year-old Kurt G\"odel
revolutionized mathematics. But even more impressive is the fact that in 1959, 25-year-old Tony Hoare
invented {\it Quicksort}, one of the most useful algorithms of all time.
Indeed, if you google "Top Ten Algorithms", you would get that Quicksort is ranked seventh.

{\bf Sorting: The stupid way}

Suppose that you are given a list of numbers, say
$$
L=[3,1,4,1,5,9,2,3] \quad,
$$
and you want to sort them in {\bf increasing} order, getting the sorted list $S$. 
The {\bf stupid} way is to first find the smallest entry, put it as the first entry of $S$,
and delete it from $L$, getting a shorter $L$. Keep doing it until $L$ is empty.

Here we used the {\bf subroutine} $min(L)$, so we need to describe it.
We also need to keep track of the location of the minimum.

The way to do it is to initialize `champ' to be the first entry, and then, going through the entries of $L$, 
{\bf compare} each of them to the current `champ', if it is smaller, you update `champ' to be that entry.

Let's find $min(L)$ for the above $L$.

$champ:=3$ . Is $1<champ$? , yes!, so the new champ is $1$, at $i=2$. Is $4<champ$?, no!, keep going. Is $1<champ$?, no, keep going.
etc. Altogether we need $n-1$ {\bf comparisons} to find the first entry of $S$. We now kick the second entry, $1$, from $L$
getting a new list with $n-1$ entries.

Iterating, we need $(n-1)+(n-2)+ \dots +1=n(n-1)/2$ comparisons, a {\bf quadratic} running time, that computer scientists
denote by $O(n^2)$.

Quicksort can do it much faster, {\it on average}. Let's recall this famous algorithm.

{\bf Input}:  A list $L$ of length $n$

{\bf Output}: The sorted list $S=Q(L)$.

\vfill\eject

$\bullet$ If the length of $L$ is $1$ or $0$ then Return $L$.

$\bullet$ Otherwise, pick $1 \leq i \leq n$ randomly, and let $p:=L[i]$.

$\bullet$ Create two shorter lists $L_1$ and $L_2$. $L_1$ consists of the entries that are  $\leq p$ and
$L_2$ consists of those that are $>p$. The output is
$$
Q(L)=[Q(L_1),p, Q(L_2)] \quad .
$$

Note that forming the two smaller lists $L_1$ and $L_2$ requires $n-1$ comparisons. If you are extremely unlucky,
one of the shorter lists will be empty, and if the unlucky streak persists you may have to do, just as with
the stupid way, $(n-1)+(n-2)+ \dots +1 \, = \, n(n-1)/2$ comparisons. Note that the
probability of that happening is pretty small: $1/n!$.

On the other hand, if you are extremely lucky, $L_1$ and $L_2$ would be of the same size $(n-1)/2$. If this
is true also for the subsequent breakings, and their descendants, then you would
get lists of length $1$ after $\log_2 \, n$ iterations. So in the {\bf best case} scenario the running time
is $n \cdot \log_2 n \, = \, O(n\, \log n )$. Unfortunately, the probability of that is also extremely small.

What about the {\bf expected} number of comparisons, if you input a random list of length $n$? 
Is the average time it takes to perform half-way between the best case scenario of $n \cdot \log_2 \, n$ and
the worst case of $n(n-1)/2$?, which is about $\frac{1}{4}\,n^2$? Thank God, not. It is  much closer to the best-possible case, and has
the same order of magnitude, namely $O(n\,\log n)$.

But let us first test Quicksort by {\it simulation}.

It takes a second to code Quicksort in Maple. Here it is:

{\tt
Q:=proc(L) local n,i,p,L1,L2,j:n:=nops(L): if n=0 or n=1 then RETURN(L): fi:
i:=rand(1..n)(): p:=L[i]: L1:=[]: L2:=[]: for j from 1 to n do  if j<>i then  if L[j]<p then   L1:=[op(L1),L[j]]:  else   L2:=[op(L2),L[j]]: fi: fi:od:
[op(Q(L1)),p,op(Q(L2))]:
end:
}

[If you have Maple, please copy-and-paste this Maple code into a Maple session, and experiment yourself!]

Using {\tt randperm} in the package {\tt combinat}, let's find the running time on $100$ random permutations of length $1000$, 
by typing

{\tt with(combinat): L:=[seq(time(Q(randperm(1000))),i=1..100)];}

The output that we got (of course, being random, you would get something different) is this:

{\tt
[0.104, 0.080, 0.084, 0.080, 0.068, 0.088, 0.076, 0.084, 0.088, 0.076, 0.080, 0.076, 0.076, 0.080, 0.084, 0.080, 0.084, 0.084, 0.076, 0.080, 
0.076, 0.076, 0.080,    0.084, 0.076, 0.072, 0.084, 0.072, 0.088, 0.080, 0.080, 0.080, 0.076, 0.080, 0.080, 0.076, 0.072, 0.080, 0.080, 0.072, 0.076, 
0.088, 0.076, 0.096, 0.076, 0.076,    0.080, 0.076, 0.080, 0.088, 0.084, 0.080, 0.084, 0.076, 0.088, 0.076, 0.076, 0.080, 0.080, 0.084, 0.076, 0.076, 
0.088, 0.076, 0.076, 0.080, 0.076, 0.084, 0.076,0.088, 0.084, 0.080, 0.072, 0.088, 0.084, 0.076, 0.072, 0.080, 0.080, 0.088, 0.084, 0.088, 
0.084, 0.076, 0.072, 0.080, 0.092, 0.080, 0.076, 0.084, 0.080, 0.076,    0.076, 0.092, 0.076, 0.080, 0.076, 0.076, 0.088, 0.080] .}

The smallest running time , $min(L)$, is $0.068$ seconds,  while the largest is $0.104$ seconds. The average, gotten by
typing {\tt convert(L,`+`)/nops(L);} is (for this random run),  $0.08028$.

Doing it with another run of $1000$ permutations (do it!)
we got that the minimum was $0.068$, the maximum was $0.100$ and the average was $0.08084$. What did you get?

{\bf The  closed-form expression for the average running time of Quicksort}

Let  $c_n$ be the expected number of comparisons performed in executing Quicksort on a random list of length $n$.

The traditional approach for proving that $c_n \, = \, O(n \, \log n)$ is presented in many textbooks 
including the  two {\bf bibles}
(both the {\it old testament} [Knu] (p.121), and the {\it new testament}, [GKP], (pp. 27-29)).
It can also be found in Wikipedia.

A nice and lucid account of the human approach to deriving an explicit expression for $c_n$ can be found in the
``The Concrete Tetrahedron'' [KaP], where it is used as a {\it motivating} example.  
They spend quite a few pages, using {\it human ingenuity}, while our derivation, to be presented in the next section,
takes a few seconds. The advantage of our approach, besides being quicker and less painful for humans,
is that it extends to the derivation of explicit expressions for higher moments, that we will describe later in this
article, where no human (without computer) has a chance.

But the {\it starting point}, setting a {\it recurrence} for $c_n$ is the same as the textbook approach,
so let us borrow this part from p. 4 of [KaP].

{\it ``In the general case, when we are sorting $n$ numbers and choose a pivot $p$, that pivot can be the
$k$-th smallest element of the list for any $k=1,\dots,n$. In any case, we need $n-1$ comparisons to bring the
$k-1$ smaller elements to the left and the $n-k$ greater elements to the right. Then we
need $c_{k-1}$ comparisons on average to sort the left part and $c_{n-k}$ comparisons on average to soft the right part,
thus $n-1+c_{k-1}+c_{n-k}$ in total. Taking the average over all possible choices for $k$, we find}
$$
c_n \, = \,
\frac{1}{n} \sum_{k=1}^{n} ((n-1)+ c_{k-1}+ c_{n-k} ) = (n-1) + \frac{1}{n} \sum_{k=1}^{n} (c_{k-1}+c_{n-k}) \quad .
= (n-1) + \frac{2}{n} \sum_{k=1}^{n} \, c_{k-1} \quad .
$$

This recurrence can be used to easily compute the first $100$ terms of the sequence of rational numbers $c_n$, and 
Fig. 1.3 of [KaP] clearly shows that $c_n$ grows much slower then $n(n-1)/2$. They conclude that section with the
remark : {\bf `` but a picture is not a proof.''} .

Since, according to them, `a picture is not a proof', Kauers and Paule [KaP], (and [GKP], and countless other
textbooks) spend a few more pages, by using {\it human-generated} {\it manipulatorics},
to rigorously derive the following {\bf closed form} formula for $c_n$.

{\bf Theorem 1} ([KaP], p.8, end of section 1.3; [GKP], Eq. (2.14), p. 29, and other places):
$$
c_n = 2(n+1) H_n - 4n \quad (n \geq 1) \quad .
$$
Here $H_n$ are the {\it Harmonic numbers}
$$
H_n:=\sum_{i=1}^{n} \frac{1}{n} \quad .
$$

Since $H_n=\log n \, + \, \gamma \, + \, o(1)$, where $\gamma= 0.5772156649\dots$ (thanks to Leonhard Euler),
it follows that indeed  $c_n=O(n\, log(n))$, and more precisely, $c_n=2n\log n + O(n)$, only 
$2/\log(2)=2.88539008\dots$  times the running time in the best-case scenario.

We will now describe our way of using the recurrence $c_n=(n-1) + \frac{2}{n} \sum_{k=1}^{n} \, c_{k-1}$ to derive
Theorem 1.

{\bf Our Derivation of Theorem 1: Keep it Simple Stupid}

What if you are not as smart and/or knowledgeable as the authors of [KaP] and [GKP]? 
And even if you are, don't you have better things to do?
Here is our `dumb' way, that may be considered as a form of `{\it machine learning}'.

First, you write a {\it one-line} Maple procedure

{\tt c:=proc(n) option remember: if n=0 then 0 else (n-1)+2/n*add(c(k),k=0..n-1) fi:end:}

[Once again, if you have Maple, please copy-and-paste this Maple code (and what follows) to a Maple session, so that you can
redo what we are doing.]

Now you make an {\it educated guess} that $c(n)$ is a polynomial of degree $1$ in both $n$ and  $H_n$, setting the
{\it template}
$$
C:= \, a \,+ \, b \cdot n \, + \, c \cdot H_n \,+ \, d \cdot nH_n \quad,
$$
with the {\it undetermined coefficients} $a,b,c,d$. Define in Maple

{\tt C:=a+b*n+c*Hn+d*n*Hn;} 

Using the first six values of $c(n)$ as the {\it training data set}, we type:

{\tt C:=subs(solve($\{$ seq(subs({n=i,Hn=add(1/j,j=1..i)},C)-c(i),i=1..6)$\}$,$\{$a,b,c,d$\}$),C);} 

and {\it lo and behold}, we get (in one nano-second) that it seems that

$$
c_n\, = \, -4\,n \,+ \, 2\,(n+1)\, H_n  \quad .
$$

But so far we only know that it is true for six cases ($1 \leq n \leq 6$). We can easily test it,
using the next $300$ cases as {\it testing data set}, and confirm that it keeps holding up to $n=306$.

Indeed, if you copy-and-paste the next line into a Maple session 

$\{$ seq(subs(n=i,Hn=add(1/j,j=1..i),C)-c(i),i=7..306)$\}$; 

You would immediately get

$$
\left\{ 0 \right\}  \quad .
$$

If you are more patient, and skeptical, you can redo it with $306$ replaced by higher numbers, but as we will soon see, that
would be wasting your computer's time.

The purist would now say: `we need to know it for {\it all} $n$, not just the first $306$ cases!'?
To him we reply, that, since, by definition, $H_n$ is defined by the recurrence
$$
H_n-H_{n-1}= \frac{1}{n} \quad , \quad H_0=0 \quad ,
$$
plugging-in into the defining recurrence for $c_n$ one would get a summation identity featuring Harmonic
numbers, handled so well by Carsten Schneider's powerful Mathematica package Sigma [S1] (see also [S2]).
Of course, in this simple case, this can be easily done by hand, but in the more complicated cases
that we will soon encounter, the fact that such a package exists gives us the peace of mind that we need.

In fact, since it is a finite calculation, it should be easy to come up with an {\it a priori} $N_0$ 
(in fact, in this case $N_0=6$ suffices), for which checking
it for the first $N_0$ cases would rigorously imply its truth for `all' $n$. 
Since it is possible to find such an $N_0$, and we are sure that it is much smaller than $306$,
why bother?

{\bf What about the Variance?}

As we all know (for example, the St. Petersburg paradox), the expectation of a random variable, while definitely
the most important number associated with a random variable, does not tell us everything about it. 
The next-in-line, in importance, is the {\it variance},
or equivalently, its square-root, called its {\it standard deviation}.

The next theorem, that we rediscovered from scratch, is the answer to Exercise 8(b) in section
6.2.2 of the Knuth's ACPIII ([Knu], pp. 448 (question); p. 672 (answer, but no proof, or even reference)). 
Knuth refers to a paper of P.F. Widley, Comp. J. {\bf 3}
(1960), 86, and mentions that Widley found a recurrence for the numerical computation of the variance, but he did not obtain
a solution). This formula also appears in [KneS], Eq. (32).

{\bf Theorem 2} (Knuth, [Knu], answer to Ex. 8(b) in section  6.2.2)): 
The variance of the random variable ``number of comparisons in Quicksort applied to lists of length $n$" is 
$$
n ( 7\,n+13 )  \, - \, 2\,(n+1)\, H_{{1}} ( n )  -4\, ( n+1 ) ^{2}H_{{2}} ( n )  \quad,
$$
where
$$
H_1(n):= \sum_{i=1}^{n} \, \frac{1}{i} \quad , \quad
H_2(n):= \sum_{i=1}^{n} \, \frac{1}{i^2} \quad .
$$
Its asymptotic expression is
$$
( 7 \,- \,\frac{2}{3}\,{\pi }^{2} ) {n}^{2}+ ( 13-2\,\ln  ( n
 ) -2\,\gamma-4/3\,{\pi }^{2} ) n-2\,\ln  ( n ) 
-2\,\gamma-2/3\,{\pi }^{2} \, + \, o(1) \quad .       
$$

As noted by Knuth ([Knu], bottom of p. 121), it follows that the distribution is {\it concentrated around the mean}. Indeed the 
asymptotic {\it coefficient of variation} is $o(1)$, but the very weak $O(1/log(n))$. 

{\bf What about the higher moments?}

As far we know, no one bothered, so far, to find explicit expressions for higher moments. 
The leading asymptotics for the third moment is given by Cramer (Eq. (2.9)), explicitly as $((16 \zeta(3)-19)+o(1))n^3$, 
and numerically for the fourth moment (Eq. (2.10))
where it is stated that is it (to eight decimal figures) $n^4\,(0.73794549+o(1))$. We found (see below) that the {\it exact value} of the leading coefficient
(of $n^4$) is
$$
\frac{4}{15} \pi^4 \, - \, 28 \pi^2 \, +  \, \frac{2260}{9} \quad .
$$

Let's define $H_m(n)$ to be the $n$-th partial sum of $\zeta(m)$:
$$
H_m(n):= \sum_{i=1}^{n} \, \frac{1}{i^m} \quad .
$$

{\bf Theorem 3}:  The third moment (about the mean) of the random variable ``number of comparisons in Quicksort applied to lists of length $n$" is 
$$
-n ( 19\,{n}^{2}+81\,n+104 ) +H_{{1}} ( n ) 
 ( 14\,n+14 ) +12\, ( n+1 ) ^{2}H_{{2}} ( n
 ) +16\, ( n+1 ) ^{3}H_{{3}} ( n ) \quad .
$$
It is asymptotic  to
$$
( -19+16\,\zeta  ( 3 )  ) {n}^{3}+ ( -81+2
\,{\pi }^{2}+48\,\zeta  ( 3 )  ) {n}^{2}+ ( -104+
14\,\ln  ( n ) 
+14\,\gamma+4\,{\pi }^{2}+48\,\zeta  ( 3
 )  ) n 
$$
$$
+14\,\ln  ( n ) +14\,\gamma+2\,{\pi }^{2}
+16\,\zeta  ( 3 ) \, + \, o(1) \quad .
$$
It follows that the limit of the scaled third moment (skewness) converges to
$$
{\frac{-19+16\,\zeta  ( 3 ) }{ ( 7-2/3\,{\pi }^{2} ) ^{3/2}}} \, = \, 0.8548818671325885\dots \quad .
$$

{\bf Theorem 4}:  The fourth moment (about the mean) of the random variable ``number of comparisons in Quicksort applied to lists of length $n$" is 
$$
1/9\,n ( 2260\,{n}^{3}+9658\,{n}^{2}+15497\,n+11357 ) -2\,
 ( n+1 )  ( 42\,{n}^{2}+78\,n+77 ) H_{{1}}
 ( n ) 
$$
$$
+12\, ( n+1 ) ^{2} ( H_{{1}} ( n
 )  ) ^{2}+ ( -4\, ( 42\,{n}^{2}+78\,n+31
 )  ( n+1 ) ^{2}+48\, ( n+1 ) ^{3}H_{{1}}
 ( n )  ) H_{{2}} ( n ) 
$$
$$
+48\, ( n+1
 ) ^{4} ( H_{{2}} ( n )  ) ^{2}-96\,
 ( n+1 ) ^{3}H_{{3}} ( n ) -96\, ( n+1
 ) ^{4}H_{{4}} ( n ) \quad .
$$
It is asymptotic  to
$$
 ( {\frac {2260}{9}}-28\,{\pi }^{2}+{\frac {4}{15}}\,{\pi }^{4}
 ) {n}^{4}+ ( {\frac {9658}{9}}-84\,\ln  ( n ) -
84\,\gamma+1/6\, ( -648+48\,\ln  ( n ) +48\,\gamma
 ) {\pi }^{2}+{\frac {16}{15}}\,{\pi }^{4}-96\,\zeta  ( 3
 )  ) {n}^{3}
$$
$$
+ ( {\frac {15497}{9}}-240\,\ln  ( n
 ) -240\,\gamma+12\, ( \ln  ( n ) +\gamma
 ) ^{2}+1/6\, ( -916+144\,\ln  ( n ) +144\,\gamma
 ) {\pi }^{2}+8/5\,{\pi }^{4}-288\,\zeta  ( 3 ) 
 ) {n}^{2}
$$
$$
+ ( {\frac {11357}{9}}-310\,\ln  ( n ) 
-310\,\gamma+24\, ( \ln  ( n ) +\gamma ) ^{2}+1/6
\, ( -560+144\,\ln  ( n ) +144\,\gamma ) {\pi }^{
2}+{\frac {16}{15}}\,{\pi }^{4}-288\,\zeta  ( 3 )  ) n
$$
$$
-154\,\ln  ( n ) -154\,\gamma+12\, ( \ln  ( n
 ) +\gamma ) ^{2}+1/6\, ( -124+48\,\ln  ( n
 ) +48\,\gamma ) {\pi }^{2}+{\frac {4}{15}}\,{\pi }^{4}-96
\,\zeta  ( 3 ) \, + \, o(1) \quad .
$$
It follows that the limit of the scaled fourth moment (kurtosis) converges to
$$
{\frac{{\frac{2260}{9}}-28\,{\pi }^{2}+{\frac{4}{15}}\,{\pi }^{4}}{
 ( 7-2/3\,{\pi }^{2} ) ^{2}}}
\, = \,
4.1781156382698542\dots \quad .
$$

{\bf Theorem 5}:  The fifth moment (about the mean) 
of the random variable ``number of comparisons in Quicksort applied to lists of length $n$" is 
$$
-{\frac {1}{108}}\,n ( 229621\,{n}^{4}+1422035\,{n}^{3}+3401325\,{
n}^{2}+3915865\,n+2217794 ) 
$$
$$
+2\, ( n+1 )  ( 190\,
{n}^{3}+1300\,{n}^{2}+1950\,n+1171 ) H_{{1}} ( n ) -
280\, ( n+1 ) ^{2} ( H_{{1}} ( n )  ) 
^{2}+ 
$$
$$
( 20\, ( 38\,{n}^{3}+204\,{n}^{2}+286\,n+91 ) 
 ( n+1 ) ^{2}-800\, ( n+1 ) ^{3}H_{{1}} ( n
 )  ) H_{{2}} ( n ) -480\, ( n+1 ) ^{
4} ( H_{{2}} ( n )  ) ^{2}
$$
$$
+ ( 80\, ( 14
\,{n}^{2}+26\,n+17 )  ( n+1 ) ^{3}-320\, ( n+1
 ) ^{4}H_{{1}} ( n ) -640\, ( n+1 ) ^{5}H_{
{2}} ( n )  ) H_{{3}} ( n ) 
$$
$$
+960\, ( n
+1 ) ^{4}H_{{4}} ( n ) +768\, ( n+1 ) ^{5}H
_{{5}} ( n ) \quad .
$$
It is asymptotic  to
$$
 ( -{\frac {229621}{108}}+{\frac {380}{3}}\,{\pi }^{2}+ ( 
1120-{\frac {320}{3}}\,{\pi }^{2} ) \zeta  ( 3 ) +768
\,\zeta  ( 5 )  ) {n}^{5}
$$
$$
+ ( -{\frac {1422035}{
108}}+380\,\ln  ( n ) +380\,\gamma+{\frac {2800}{3}}\,{\pi }
^{2}-8/3\,{\pi }^{4}+ ( 5440-320\,\ln  ( n ) -320\,
\gamma-{\frac {1600}{3}}\,{\pi }^{2} ) \zeta  ( 3 ) +
3840\,\zeta  ( 5 )  ) {n}^{4}
$$
$$
+ ( -{\frac {125975}
{4}}+2980\,\ln  ( n ) +2980\,\gamma+1/6\, ( 14640-800\,
\ln  ( n ) -800\,\gamma ) {\pi }^{2}-{\frac {32}{3}}\,
{\pi }^{4}+ 
$$
$$
( 10960-1280\,\ln  ( n ) -1280\,\gamma-{
\frac {3200}{3}}\,{\pi }^{2} ) \zeta  ( 3 ) +7680\,
\zeta  ( 5 )  ) {n}^{3}
$$
$$
+ ( -{\frac {3915865}{108}
}+6500\,\ln  ( n ) +6500\,\gamma-280\, ( \ln  ( n
 ) +\gamma ) ^{2}+1/6\, ( 17340-2400\,\ln  ( n
 ) -2400\,\gamma ) {\pi }^{2}-16\,{\pi }^{4}
$$
$$
+ ( 11440-
1920\,\ln  ( n ) -1920\,\gamma-{\frac {3200}{3}}\,{\pi }^{2}
 ) \zeta  ( 3 ) +7680\,\zeta  ( 5 ) 
 ) {n}^{2}
$$
$$
+ ( -{\frac {1108897}{54}}+6242\,\ln  ( n
 ) +6242\,\gamma-560\, ( \ln  ( n ) +\gamma
 ) ^{2}+1/6\, ( 9360-2400\,\ln  ( n ) -2400\,
\gamma ) {\pi }^{2}-{\frac {32}{3}}\,{\pi }^{4}
$$
$$
+ ( 6160-1280
\,\ln  ( n ) -1280\,\gamma-{\frac {1600}{3}}\,{\pi }^{2}
 ) \zeta  ( 3 ) 
$$
$$
+3840\,\zeta  ( 5 ) 
 ) n+2342\,\ln  ( n ) +2342\,\gamma-280\, ( \ln 
 ( n ) +\gamma ) ^{2}+1/6\, ( 1820-800\,\ln 
 ( n ) 
$$
$$
-800\,\gamma ) {\pi }^{2}-8/3\,{\pi }^{4}+
 ( 1360-320\,\ln  ( n ) 
$$
$$
-320\,\gamma-{\frac {320}{3}}\,
{\pi }^{2} ) \zeta  ( 3 ) +768\,\zeta  ( 5
 ) \, + \, o(1) \quad . 
$$

It follows that the limit of the scaled fifth moment  converges to
$$
{\frac{-{\frac{229621}{108}}+{\frac{380}{3}}\,{\pi }^{2}+ \left( 
1120-{\frac{320}{3}}\,{\pi }^{2} \right) \zeta  \left( 3 \right) +768
\,\zeta  \left( 5 \right) }{ \left( 7-2/3\,{\pi }^{2} \right)^{5/2}}}
\, = \, 10.64616337467387850\dots \quad .
$$

{\bf Theorem 6}:  The sixth moment (about the mean) 
of the random variable ``number of comparisons in Quicksort applied to lists of length $n$" is 
$$
{\frac{1}{2700}}\,n ( 74250517\,{n}^{5}+523547007\,{n}^{4}+
1579578725\,{n}^{3}+2571768745\,{n}^{2}+2342670258\,n+1133389148
 ) 
$$
$$
-2/3\, ( n+1 )  ( 11300\,{n}^{4}+56270\,{n}^{3
}+135760\,{n}^{2}+145510\,n+68427 ) H_{{1}} ( n ) +20
\, ( 63\,{n}^{2}+117\,n+329 )  ( n+1 ) ^{2}
 ( H_{{1}} ( n )  ) ^{2}
$$
$$
-120\, ( n+1
 ) ^{3} ( H_{{1}} ( n )  ) ^{3}
$$
$$
+ ( -4/
3\, ( 11300\,{n}^{4}+51710\,{n}^{3}+101830\,{n}^{2}+93640\,n+26013
 )  ( n+1 ) ^{2}+240\, ( 21\,{n}^{2}+39\,n+68
 )  ( n+1 ) ^{3}H_{{1}} ( n ) 
$$
$$
-720\,
 ( n+1 ) ^{4} ( H_{{1}} ( n )  ) ^{2}
 ) H_{{2}} ( n ) 
$$
$$
+ ( 240\, ( 21\,{n}^{2}+39
\,n+37 )  ( n+1 ) ^{4}-1440\, ( n+1 ) ^{5}H
_{{1}} ( n )  )  ( H_{{2}} ( n ) ) ^{2}
$$
$$
-960\, ( n+1 ) ^{6} ( H_{{2}} ( n
 )  ) ^{3} 
$$
$$
+ ( -160\, ( 38\,{n}^{3}+225\,{n}^{2}+325\,n+159 )  ( n+1 ) ^{3}+7360\, ( n+1 )^{4}H_{{1}} ( n ) 
$$
$$
+9600\, ( n+1 ) ^{5}H_{{2}}
 ( n )  ) H_{{3}} ( n ) +2560\, ( n+1
 ) ^{6} ( H_{{3}} ( n )  ) ^{2}+ ( -
480\, ( 21\,{n}^{2}+39\,n+37 )  ( n+1 ) ^{4}+2880
\, ( n+1 ) ^{5}H_{{1}} ( n ) 
$$
$$
+5760\, ( n+1
 ) ^{6}H_{{2}} ( n )  ) H_{{4}} ( n
 ) -11520\, ( n+1 ) ^{5}H_{{5}} ( n ) -7680
\, ( n+1 ) ^{6}H_{{6}} ( n )  \quad .
$$
Its asymptotic  expressions can be found here:

{\tt http://sites.math.rutgers.edu/\~{}zeilberg/tokhniot/oQuickSortAnalysis3.txt} \quad .

It follows that the limit of the scaled sixth moment  converges to
$$
{\frac {{\frac {74250517}{2700}}-{\frac {22600}{9}}\,{\pi }^{2}+140\,{
\pi }^{4}-{\frac {88}{7}}\,{\pi }^{6}-6080\,\zeta  \left( 3 \right) +
2560\, \left( \zeta  \left( 3 \right)  \right) ^{2}+{\frac {1}{90}}\,
 \left( 960\,{\pi }^{2}-10080 \right) {\pi }^{4}}{ \left( 7-2/3\,{\pi }
^{2} \right) ^{3}}}
$$
$$
\, = \, 44.42707770816977761\dots \quad .
$$

{\bf Theorem 7}:  The seventh moment (about the mean) 
of the random variable ``number of comparisons in Quicksort applied to lists of length $n$", as well as its asymptotics, can be found here:

{\tt http://sites.math.rutgers.edu/\~{}zeilberg/tokhniot/oQuickSortAnalysis3.txt} \quad .

It follows that the limit of the scaled seventh moment  converges to
$$
\frac{1}{81000} \cdot (7-2/3\,{\pi }^{2})^{-\frac{7}{2}} \cdot
$$
$$
(-30532750703+2411020500\,{\pi }^{2}-14364000\,{\pi }^{4}+11390400000\,\zeta  \left( 3 \right) -1270080000\,\zeta  \left( 3 \right) {\pi }^{2}+12096000\,
\zeta  \left( 3 \right) {\pi }^{4}
$$
$$
-870912000\,\zeta  \left( 5 \right) {\pi }^{2}+9144576000\,\zeta  \left( 5 \right) +7464960000\,\zeta  \left( 7 \right) 
)
$$
$$
\, = \, 
 179.7219197356178684\dots \quad .
$$

{\bf Theorem 8}:  The eighth moment (about the mean) 
of the random variable ``number of comparisons in Quicksort applied to lists of length $n$", as well as its asymptotics, can be found here:

{\tt http://sites.math.rutgers.edu/\~{}zeilberg/tokhniot/oQuickSortAnalysis3.txt} \quad .

It follows that the limit of the scaled eighth moment  converges to (in Maple format)
$$
\frac{1}{183750} \cdot(7-2/3\,{\pi }^{2})^{-4} \cdot
$$
$$
(\, 90558126238639-7640378199300\,{\pi }^{2}+69766200000\,{\pi }^{4}-5556600000\,{\pi }^{6}-354564000\,{\pi }^{8}-28353601080000\,\zeta  ( 3 ) +
$$
$$
1689206400000\,\zeta  ( 3 ) {\pi }^{2}+7468070400000\, ( \zeta  ( 3 )  ) ^{2}-711244800000\, ( \zeta  ( 3
 )  ) ^{2}{\pi }^{2}
$$
$$
-12162286080000\,\zeta  ( 5 ) +10241925120000\,\zeta  ( 5 ) \zeta  ( 3 ) \, )
$$
$$
= \,858.203203990002260\dots \quad .
$$

{\bf How Theorems 2-8 were discovered}

The same way as we discovered Theorem 1! Even more informative than the moments is the full discrete probability
distribution of the random variable `number of comparisons in Quicksort applied to permutations of length $n$'.
Let's call it $X_n$.

It is useful to introduce the {\bf probability generating function}
$$
g_n(t) \, := \, \sum_{k=0}^{n(n-1)/2} Pr(X_n=k)\,t^k \quad .
$$
The same reasoning that lead to the recurrence for $c_n$ yields the recurrence (ans. to ex. 8(a) of section 6.2.2. of [Knu])
$$
g_n(t) \, = \, \frac{t^{n-1}}{n} \sum_{k=1}^{n} g_{k-1}(t) g_{n-k}(t) \quad .
\eqno(DEK)
$$

Recall that the $r$-th moment is given in terms of the probability generating function
$$
E[X_n^r]= (t\frac{d}{dt})^r g_n(t) \, \vert_{t=1} \quad .
$$
More informative is the moment-about-the-mean
$$
m_r(n) := E[(X_n - c_n)^r] \quad ,
$$
that can be easily derived from the straight moments $\{ E[X_n^l] \,|\, 1 \leq l \leq r\}$, using the Binomial theorem and
linearity of expectation.

Our {\it data driven} approach is to use Maple to generate as many terms of the sequence of polynomials $g_n(t)$ that
it would care to give us. See

{\tt http://sites.math.rutgers.edu/\~{}zeilberg/tokhniot/oQuickSortAnalysis1.txt}

for the first $130$ of them. As you can see, they get pretty large, and it is important to keep them in rational arithmetic.

Once you have them, it is very easily, for any desired moment, $m_r(n)$, to get the first $130$ (or whatever) terms
of the {\bf numerical} sequence $\{m_r(n)\}$. Now you make the {\it educated guess} that there exists
a polynomial expression for $m_r(n)$ in terms of $n$ and
$$
H_m(n) := \sum_{i=1}^{n} \frac{1}{i^m} \quad ,
$$
for $1 \leq m \leq r$. In other words, there exists a polynomial of $r+1$ variables, let's call it $F_r(x_0, x_1, \dots ,x_r)$ such that
$$
m_r(n) \, = \, F_r( n, H_1(n), \dots, H_r(n) ) \quad .
$$
We first try a polynomial of (total) degree $1$, then $2$, until we succeed. Just like in the case for $c_n$, we write 
a generic $F_r$, of the specified degree in terms of its {\it undetermined coefficients}, plug-in enough data to get
a few more equations than unknowns, and then solve them. Once we have a conjecture, we test it for quite a few more data
points, getting a very plausible conjecture.

How do we prove them? The non-linear
recurrence for the probability generating function $g_n(t)$, implies extremely complicated recurrences for
the moments, where the recurrence for a specific moment $m_r(n)$ involves lower moments $m_s(n)$ ($1 \leq s <r$) that we already know.
There are {\it decidable}, using (for example) Carsten Schneider's Mathematica package [S1].
Since we have the option to have it proved rigorously, why bother? With all due respect to Theorem 8 above,
it is not important enough to have a fully rigorous proof. A semi-rigorous proof obtained by checking
sufficiently many special cases is good enough for us.

{\bf Getting more data for the moments}

If we are only interested in the first few moments, say, the first $20$, then it is wasteful to compute the full $g_n(t)$.
We write  $t=1+w$ and use the fact that
$$
g_n(1+w) \, = \, \sum_{r=0}^{\infty} \frac{f_r(n)}{r!} w^r \quad ,
$$
where $f_r(n)$ are the {\bf factorial moments}, from which the straight moments $E[X_n^r]$, and hence the moments-about-the-mean,
$m_r(n)$, can be computed. 

The non-linear recurrence $(DEK)$ implies that
$$
g_n(1+w) \, = \, \frac{(1+w)^{n-1}}{n} \sum_{k=1}^{n} g_{k-1}(1+w) g_{n-k}(1+w) \quad .
\eqno(DEK')
$$
If we are only interested in the first $M$ factorial moments, we can truncate at each step and only keep the first $M$ coefficients in $w$,
and get much further.

{\bf Approximating the Limiting Distribution using Symbol-Crunching}

Since we have a closed form expression for both the expectation, $c_n$, and the variance $m_2(n)$, we can form the
scaled, distribution
$$
Z_n \, := \, \frac{X_n -c_n}{\sqrt{m_2(n)}} \quad ,
$$
and ask about the {\it limiting distribution} $\lim_{n \rightarrow \infty} Z_n$. Of course its expectation is $0$ and its
variances is $1$, and the exact values of its $3$-rd through $8$-th moments were given in Theorems 3 through to 8, but here there are again,
in floating point approximations:
$$
[.85488186713258853660, 4.1781156382698542397, 10.646163374673878503, 44.427077708169777614, 
$$
$$
179.72191973561786840, 858.20320399000226017] \quad .
$$

Using symbolic computation, one can get a very good approximation for the limiting distribution $Z_{\infty}$ by taking,
say $Z_{130}$. See the diagram in

{\tt http://sites.math.rutgers.edu/\~{}zeilberg/tokhniot/qsort/pdf130.html} 

for the {\bf density function}.

Since we know the exact values of $c_n$ and $m_2(n)$, we can get a very good approximation of $X_n$ for very large $n$
by using the approximation $c_n + \sqrt{m_2(n)}Z_{130}$ instead of  $c_n + \sqrt{m_2(n)}Z_{\infty}$, and compute very
good approximation for the probability that quicksort will take more than a specified number of comparisons. See an
example for $n=10000$ in 

{\tt http://sites.math.rutgers.edu/\~{}zeilberg/tokhniot/qsort/cmf10000.html}

{\bf The Maple package QuickSortAnalysis.txt}

Everything in this paper was done by the first author by running the Maple package \hfill\break
{\tt QuickSortAnalysis.txt}
written by the second author. It is available from the front of this article

{\tt http://sites.math.rutgers.edu/\~zeilberg/mamarim/mamarimhtml/qsort.html} 

{\bf Conclusion}

As with most of our joint papers, more important than the actual results is the illustration of a {\bf methodology}
of experimental mathematics, based on `guessing' and `big data', that forms an alternative to traditional
human-generated analysis of algorithms using ad-hoc manipulations. It also forms an alternative to mere
simulations, by getting exact expressions for the moments, and good approximations for the probability distributions.

{\bf References}

[C] Michael Cramer, {\it A note concerning the limit distribution of the quicksort algorithm},
Informatique The\'eriques et Applications, {\bf 30} (1996), 195-207.

[GKP] Ronald L. Graham, Donald E. Knuth, and Oren Patashnik, {\it ``Concrete Mathematics''},  Addison-Wesley, 1989.

[KaP] Manuel Kauers  and Peter Paule, {\it ``The Concrete Tetrahedron''}, Springer, 2011.

[KneS] Charles Knessl and Wojciech Szpankowski, {\it Quicksort algorithm again revisited},
Discrete Mathematics and Theoretical Computer Science, {\bf 3} (1999), 43-64.

[Knu] Donald E. Knuth, {``The Art of Computer Programming'', Volume 3: Sorting and Searching}, Addison-Wesley, 1973.

[S1] Carsten Schneider, {\it The Summation package Sigma}, A Mathematica package available from \hfill\break
{\tt https://www3.risc.jku.at/research/combinat/software/Sigma/index.php} 

[S2] Carsten Schneider, {\it  Symbolic Summation Assists Combinatorics}, Sem.Lothar.Combin. {\bf 56}(2007),Article B56b (36 pages). \hfill\break
{\tt https://www3.risc.jku.at/research/combinat/software/Sigma/pub/SLC06.pdf}  

\vfill\eject

\bigskip
\hrule
\bigskip
Shalosh B. Ekhad, c/o D. Zeilberger, Department of Mathematics, Rutgers University (New Brunswick), Hill Center-Busch Campus, 110 Frelinghuysen
Rd., Piscataway, NJ 08854-8019, USA. \hfill\break
Email: {\tt ShaloshBEkhad at gmail dot com}   \quad .
\bigskip
Doron Zeilberger, Department of Mathematics, Rutgers University (New Brunswick), Hill Center-Busch Campus, 110 Frelinghuysen
Rd., Piscataway, NJ 08854-8019, USA. \hfill\break
Email: {\tt DoronZeil at gmail  dot com}   \quad .
\bigskip
\hrule
\bigskip
Exclusively published in the Personal Journal of Shalosh B.  Ekhad and Doron Zeilberger and arxiv.org \quad .
\bigskip
First Written: March 8, 2019. 
\end